\documentclass[12pt,oneside]{amsart}
\def\indicator{{\mathbf 1}}
\newtheorem{theorem}{Theorem}[section]
\usepackage[letterpaper,left=2cm,right=2cm,top=3cm,bottom=3cm]{geometry}

\newtheorem{definition}[theorem]{Definition}
\newtheorem{lemma}[theorem]{Lemma}

\renewcommand{\epsilon}{\varepsilon}
\numberwithin{equation}{section}
\begin{document}

\title{Balls are maximizers of the Riesz-type functionals with supermodular integrands}
\author{Hichem Hajaiej}

\address{Justus-Liebig-Universit\"at Giessen\\Mathematisches Institut\\Arnd Str 2, 35392 Giessen\\Germany}\email{hichem.hajaiej@gmail.com}

\begin{abstract}
For a large class of supermodular integrands, we establish conditions under which balls are the unique (up to translations) maximizers of the Riesz-type functionals with constraints. 
\end{abstract}

\maketitle
\pagestyle{headings}
\section{Introduction}
Over the last decades, one field of intense research activity has been the study of extremals of integral functionals.  The Riesz-type kind has attracted growing attention and played a crucial role in the resolution of Choquard's conjecture in a breakthrough paper by E.~H. Lieb~\cite{c1}.  The determination of cases of equality in the Riesz-rearrangement inequality has also received a large amount of interest from mathematicians due to its connection with many other functional inequalities and its several applications to physics~\cite{c2,c3,c4}.  Variational problems for steady axisymmetric vortex-rings in which kinetic energy is maximized subject to prescribed impulse involves Riesz-type functionals with constraints.  In~\cite{c5}, G.~R. Burton has proved the existence of maximizers in an extended constraint set, he has also showed that the maximizer is Schwarz symmetric (up to translations).  His method hinges on a resolution of an optimization of a Riesz-type functional under constraint~\cite[Proposition~8]{c5}.  The purpose of this paper is to answer the more general question: When do maximizers of the Riesz-type functional inherit the symmetry and monotonicity properties of the integrand involved in it? 

The method of G.~R. Burton~\cite{c5} cannot apply to solve the above problem.  In this paper, we develop a self-contained approach.  Let us give here a foretaste of our ideas.  First, we recall that: 

A Riesz-type functional is a functional of the form: 
\begin{equation*}
R(f,g)=\int_{\mathbb R^n}\int_{\mathbb R^n}\Psi\left(f(x), g(y)\right)\, r(x,y)\, dx\, dy. 
\end{equation*}

In this paper, we will consider $r(x,y)=j\left(|x-y|\right)$.  We are interested in the following maximization problem: \newline
(P1)\hfill$\sup\limits_{(f, g)\in C}J(f, g)$\hfill\mbox{~}\\
where 
\begin{equation}
J(f, g)=\int_{\mathbb R^n}\int_{\mathbb R^n}\Psi\left(f(x), g(y)\right)\, j\left(|x-y|\right)\, dx\, dy. 
\end{equation}
and 
\begin{equation}
C=(f,g): 
\begin{cases}
f:&\mathbb R^n\rightarrow\mathbb R; 0\leq f\leq k_1 \mbox{ and } \int_{\mathbb R^n}f\leq \ell_1\\
g:&\mathbb R^n\rightarrow\mathbb R; 0\leq g\leq k_2 \mbox{ and } \int_{\mathbb R^n}g\leq \ell_2
\end{cases}
\end{equation}
$\ell_1, k_1, \ell_2, k_2$ are positive numbers. 

For supermodular operators $\Psi$ and nonincreasing functions $j$, we know that $J(f, g)\leq J(f^*, g^*)$~\cite[Theorem~1]{c4}, where $u^*$ denotes the Schwarz symmetrization of $u$.  Hence the problem reduces to: \newline
(P2)\hfill$\sup\limits_{(f^*, g^*)\in C}J\left(f^*, g^*\right)$\hfill\mbox{~}\\
For continuous integrands $\Psi$ having the N-Luzin property (for any subset N having Lebesgue measure zero, $\Psi$(N) has the same property), lemma~\ref{lem23} enables us to assert that (P2) is equivalent to an optimization of a Hardy-Littlewood type functionals where balls are maximizers.  We will then extend this study to supermodular non-continuous bounded functions $\Psi$ thanks to the decomposition of these functions into $\tilde{\Psi}\left(\varphi_1(s_1), \varphi_2(s_2)\right)$ in the spirit of~\cite{c4, c6}.  The approximation of unbounded supermodular functions by bounded ones inheriting the monotonicity properties will enable us to prove that balls are maximizers in the general case. 
\\

\parindent0ex\underline{\bf Main Result:} \\
\\
Let $\Psi: \mathbb R_+\times\mathbb R_+\rightarrow\mathbb R$ be a H-Borel function satisfying: \newline
($\Psi$1)\hfill$\Psi$ vanishes at hyperplanes; \hfill\mbox{~}\\
($\Psi$2)\hfill$\Psi(b, d)-\Psi(b, c)-\Psi(a, d)+\Psi(a, c)\geq 0$ for all $0\leq a<b$ and $0\leq c <d$; \hfill\mbox{~}\\
($\Psi$3)(i)\hfill $\Psi(tx, b_2)-t\Psi(x, b_2)-\Psi(tx, b_1)+t\Psi(x, b_1)\leq 0$ for all $x\geq 0$, $0\leq b_1<b_2$ and $0<t<1$; \hfill\mbox{~}\\
($\Psi$3)(ii)\hfill$\Psi(a_2, ty)-t\Psi(a_2, y)-\Psi(a_1, ty)+t\Psi(a_1, y)\leq 0$ for all $y\geq 0$, $0\leq a_1<a_2$ and $0<t<1$; \hfill\mbox{~}\\
(j1)\hfill$j$ is nonincreasing. \hfill\mbox{~}\\

Suppose in addition that $\Psi$ is continuous with respect to each variable and has the N-Luzin property, then for all $(f_1, f_2)\in C$
\begin{equation*}
J\left(f_1, f_2\right) \leq J\left(k_1\mathrm{1}_{B_1}, k_2\mathrm{1}_{B_2}\right)
\end{equation*}
where $B_1$ and $B_2$ are centered in the origin, $\mathrm{1}_B$ is the characteristic function of $B$, and $\mu(B_1)=\ell_1/k_1$, $\mu(B_2)=\ell_2/k_2$. Moreover, if ($\Psi$2) and ($\Psi$3) hold with strict inequality, $j$ is strictly decreasing and $J\left(f_1, f_2\right)<\infty$ for any $(f_1, f_2)\in C$, then (P1) is attained by exactly two couples $\left(k_1\mathrm{1}_{B_1}, k_2\mathrm{1}_{B_2}\right)$ and $(h_1, h_2)$ where $h_1$ and $h_2$ are translates by the same vector of $k_1\mathrm{1}_{B_1}$ and $k_2\mathrm{1}_{B_2}$ (respectively). 

\section{Notations and preliminaries}
\begin{definition}
 If $A\subset\mathbb R^n$ is a measurable set of finite Lebesgue measures $\mu$, we define $A^*$, the symmetric rearrangement of the set $A$ to be the open ball centered at the origin whose volume is that of $A$, thus $A^*=\left\{x\in\mathbb R^n: |x|<r\right\}$ with $V_n r^n=\mu(A)$, $V_n$ is a constant.  
\end{definition}

For a nonnegative measurable function $u$ on $\mathbb R^n$, we require $u$ to vanish at infinity in the sense that all its positive level sets $\{x\in\mathbb R^n:u(x)>t\}$ having finite measure for $t>0$.  The set of these functions is denoted by $F_n$.  The symmetric decreasing rearrangement $u^*$ of $u$ is the unique upper semicontinuous, nonincreasing radial function that is equimeasurable with $u$.  Explicitly, $u^*(x)=\int\limits_0^\infty \indicator^*_{\{u>t\}}(x)\, dt$ where $\indicator^*_A=\indicator_{A^*}$.  We say that $u$ is Schwarz symmetric if $u\equiv u^*$. 

\begin{definition}
A reflexion $\sigma$ on $\mathbb R^n$is an isometry with the properties: 
\begin{itemize}
\item[(i)]$\sigma^2_x=\sigma_x \circ \sigma_x=x$ for all $x\in\mathbb R^n$; 
\item[(ii)]the fixed point set of $H_0$ of $\sigma$ separates $\mathbb R^n$ into two half spaces $H_+$ and $H_-$ that are interchanged by $\sigma$; 
\item[(iii)]$|x-x'|<|x-\sigma_{x'}|$ for all $x, x'\in H_+$. 
\end{itemize}
$H_+$ is the half space containing the origin. 
\end{definition}

The two point rearrangement or polarization of a real valued function $u$ with respect to a reflection $\sigma$ is defined by: 
\begin{equation}
 u^{\sigma_x}=\begin{cases}
              \max\{u(x), u(\sigma_x)\}, x \in H_+ \cup H_0,\\ 
              \min\{u(x), u(\sigma_x)\}, x \in H_-. 
             \end{cases}
\end{equation}

\begin{lemma}\label{lem21}
 Let $j:[0, \infty) \rightarrow \mathbb R$ be a nonincreasing function then $\nu(x)=\int_{{\mathbb R}^n}j\left(|x-y|\right)h(y)\, dy$ is radial and radially decreasing for any Schwarz symmetric function $h$.  If in addition $j$ is strictly radially decreasing then $\nu$ also inherits this property. 
\end{lemma}
\parindent0ex\textbf{Proof:} we will use~\cite[Lemma 2.8]{c7}: $u=u^* \Leftrightarrow u=u^\sigma$ for all $\sigma$.  It is sufficient to prove that $u(x)\geq u(\sigma_x)$ for all $x\in\mathbb R^n$, all $\sigma$.
\begin{eqnarray*}
u(x)&=&\int_{H^+}j\left(|x-y|\right)h(y)+j\left(|x-\sigma_y|\right)h(\sigma_y)\, dy\\
u(\sigma_x)&=&\int_{H^+}j\left(|\sigma_x-y|\right)h(y)+j\left(|\sigma_x-\sigma_y|\right)h(\sigma_y)\, dy\\
u(x)-u(\sigma_x)&=&\int_{H_+}j\left(|x-y|\right)[h(y)-h(\sigma_y)]-j\left(\sigma_x-y\right)[h(y)-h(\sigma_y)]\, dy\\
&=&\int_{H_+}\left(j\left(|x-y|\right)-j(\sigma_x-y)\right)\left(h(y)-h(\sigma_y)\right)\, dy
\end{eqnarray*}
By (iii) $|x-y|<|\sigma_x-y|$, it follows that $j\left(|x-y|\right)\geq j\left(|\sigma_x-y|\right)$.  On the other hand $h$ is Schwarz symmetric, hence $h(y)\geq h(\sigma_y)$ for all $y\in H_+$, the conclusion follows. 

\begin{definition}
Let $\Psi: \mathbb R_+\times\mathbb R_+ \rightarrow \mathbb R$: \\
\begin{itemize}
\item[(a)]$\Psi$ is supermodular if $(\Psi2)$ holds. 
\item[(b)]We say that $\Psi$ vanishes at hyperplanes if $\Psi(s_1, 0)=\Psi(0, s_2)=0$ for all $s_1, s_2\geq 0$. 
\end{itemize}
\end{definition}
An important property of functions satisfying (c) is that the composition
$(x, y)\mapsto \Psi\left(f(x), g(y)\right)$ is measurable on 
$\mathbb{R}_+$ for every $f, g\in F_n$.  Hence 
$j\left(|x-y|\right)\Psi\left(f(x), g(y)\right)$ is measurable 
on $\mathbb{R}_+ \times \mathbb{R}_+$. 

In the spirit of~\cite{c4} and~\cite{c6}, we obtain: 
\begin{lemma}\label{lem22}
Assume that $\Psi: \mathbb{R}_+\times \mathbb{R}_+\rightarrow\mathbb{R}$
is a supermodular bounded function vanishing at hyperplanes. 
Then there exist two bounded nondecreasing functions 
$\varphi_1$ and $\varphi_2$ on $\mathbb{R}_+$ with 
$\varphi_i(0)=0$ and a Lipschitz continuous function $\tilde\Psi$ on $\mathbb{R}_+^2$ such that $\Psi(u, v)=\tilde{\Psi}\left(\varphi(u), \varphi(v)\right)$. 
\end{lemma}
\parindent0ex\textbf{Proof:} First, we will prove the following: If $\varphi$ is a nondecreasing real-valued function defined on an interval $I$, then for every $f$ on I satisfying $|f(u)-f(v)|<c\left(\varphi(v)-\varphi(u)\right)$ where $u<v\in I$, $c$ is a constant, there exists a Lipschitz continuous function $\tilde{f}:\mathbb{R}\rightarrow[\inf{f}, \sup{f}]$ such that $f(x)=\tilde{f}\circ\varphi(x)$ (2.0). If $f$ is nondecreasing then $\tilde{f}$ is nondecreasing also. 

The result is obvious for $t=\varphi(v)$ and $s=\varphi(u)<t$ since we have $$|\tilde{f}(t)-\tilde{f}(s)|=\left|f\left(\varphi(v)\right)-f\left(\varphi(u)\right)\right|\leq c\left(\varphi(v)\right)-\left(\varphi(u)\right)=c (t-s).$$ 

Now $\tilde{f}$ has a unique extension to the closure of the image and the complement consists of a countable number of disjoint bounded intervals, it is sufficient to interpolate $\tilde{f}$ linearly between the values, that were assigned to end-points.  By construction $f=\tilde{f} \circ \varphi$ and $\tilde{f}(\mathbb R)=[\inf f, \sup f]$ the extension we have made by linear interpolation preserves of course the modulus of continuity of $\tilde{f}$: $|\tilde f(t)-\tilde f(s)|\leq c (t-s)$ for all $t>s$. If $f$ is nondecreasing, it is easy to check that this property is inherited by $\tilde f$. 

Now we can prove our lemma: 

First note that the fact that $\Psi$ is supermodular and vanishes at hyperplanes imply that 
it is nondecreasing with respect to each variable and it is nonnegative.  Now set $\varphi_1(u)=\lim\limits_{u\rightarrow +\infty}\Psi(u, v)$. 

$\varphi_1$ is well-defined on $\mathbb R_+$ since $\Psi$ is bounded and nondecreasing in the second variable.  By the supermodularity of $\Psi$, it follows that 
\begin{equation*}
\Psi(u+h_1, v+h_2)-\Psi(u, v+h_2)-\Psi(u+h_1, v)+\Psi(u, v)\geq 0
\end{equation*}
for any $u, v, h_1$ and $h_2\geq 0$. 

Letting $h_2$ tend to infinity, we obtain 
\begin{equation*}
\varphi_1(u+h_1)-\varphi(u)\geq\Psi(u+h_1, v)-\Psi(u, v)\geq 0
\end{equation*}
for all $u, v, h_1\geq 0$. 

For a fixed $v$, the last inequality enables us to apply (2.0) to $\Psi(\cdot, v)$.  Hence, there exists $\Psi^1$ such that: $\Psi(u, v)=\Psi^1(\varphi_1(u), v)$.  A moment's consideration shows that $\Psi^1$ inherits all the properties of $\Psi$.  Now set $\varphi_2(v)=\lim\limits_{u\rightarrow +\infty}\Psi(u, v)$, a similar argument ensures us that there exists $\tilde\Psi$ such that 
$\Psi^1\left(\varphi_1(u), u\right)=\tilde\Psi^1\left(\varphi_1(u), \varphi_2(v)\right)$. 

$\tilde\Psi$ has the same monotonicity property as $\Psi^1$ and consequently as $\Psi$.  Note that $\varphi_1(0)=\varphi_2(0)=0$ and the monotonicity properties of $\Psi$ imply that 
$\varphi_1$ and $\varphi_2$ are nondecreasing.  

\begin{lemma}\label{lem23}
Let $l, k>0$, $D=\{h:\mathbb R^n\rightarrow\mathbb R: 0\leq h(x)\leq k\mbox{ and }\int_{\mathbb R^n}h(x)\, dx\leq l\}$.  Suppose that $\Gamma: \mathbb R_+\rightarrow\mathbb R$ is a function satisfying:
\begin{enumerate}
 \item $\Gamma(0)=0$, 
 \item $\Gamma(t x)\leq t\Gamma(x)$ for all $x\geq 0$ and $0<t<1$. 
 \item[] Assume also that 
 \item $u: \mathbb R^n\rightarrow\mathbb R$ is a Schwarz symmetric function.  Then for every $\nu\in D: \int_{\mathbb R^n}u(x)\Gamma\left(\nu(x)\right)\, dx\leq\int_{\mathbb R^n}u(x)\Gamma\left(k \indicator_B(x)\right)\, dx$ where $B$ is the ball centered at the origin with $\mu(B)=\ell/k$. 
\end{enumerate}
\end{lemma}
\parindent0ex\textbf{Proof:} (2) implies that 
\begin{eqnarray*}
\int_{\mathbb R^n}u(x)\Gamma\left(\nu(x)\right)\, dx\leq\int_{\mathbb R^n}u(x)\Gamma(k)\frac{\nu(x)}{k}\, dx&=&
\Gamma(k)
\left[
\int_Bu(x)\left[\frac{\nu(x)}{k}-1+1\right]\, dx
+\int_{\mathbb R^n-B}\frac{u(x)\nu(x)}{k}\, dx
\right]\\
&=&\int_{\mathbb R^n}u(x)\Gamma\left(k \indicator_{B(x)}\right)\, dx\\
&&\quad+\Gamma(k)
\left[
 \int_Bu(x)\left[\frac{\nu(x)}{k}-1\right]\, dx
+\int_{\mathbb R^n-B}\frac{u(x)\nu(x)}{k}\, dx
\right]\, dx. 
\end{eqnarray*}

Using (3), it follows that the above integrals are $\leq\int_{\mathbb R^n}u(x)\Gamma\left(k \indicator_{B(x)}\right)\, dx+\Gamma(k)u(r)\left[\int_{\mathbb R^n}\frac{\nu(x)}{k}\, dx-\mu(B)\right]$ where $\mu(B)=V_rr^n$ (see definition 2.1). Thus 
$\int_{\mathbb R^n}u(x)\Gamma\left(\nu(x)\right)\, dx\leq\int_{\mathbb R^n}u(x)\Gamma\left(k \indicator_{B(x)}\right)\, dx$, since $\int_{\mathbb R^n}\frac{\nu(x)}{k}\, dx\leq \mu(B)=\ell/k$. 

If additionally $\int_{\mathbb R^n}u(x)\Gamma\left(\nu(x)\right)\, dx<\infty$ for any $\nu\in D$, (2) holds with strict inequality and $u$ is strictly decreasing, we can prove that for every $\nu\in D$: $\int_{\mathbb R^n}u(x)\Gamma\left(\nu(x)\right)\, dx<\int_{\mathbb R^n}u(x)\Gamma\left(k \indicator_{B(x)}\right)\, dx$. 

\section{Proof of the result}
For the convenience of the reader, the proof will be divided in three parts. \\
\textbf{First part:} We suppose that: $\Psi(\cdot, s_2)$ is absolutely continuous for every $s_2\geq 0$, and $\Psi(s_1, \cdot)$ is absolutely continuous for every $s_1\geq 0$.  

First note that $(\Psi1)$ and $(\Psi2)$ imply that $\Psi$ is a non-decreasing function with respect to each variable and it is nonnegative.  

Let $(f_1, f_2)\in C$, $(\Psi2)$ and (j1) imply that
\begin{eqnarray*}
J\left(f_1, f_2\right)\leq J\left(f_1^*, f_2^*\right)&=&\int_{\mathbb R^n}\int_{\mathbb R^n}\Psi\left(f_1^*(x), f_2^*(y)\right)j\left(|x-y|\right)\, dx\, dy\\
&=&\int_{\mathbb R^n}\int_{\mathbb R^n}\left(
\int_0^{f_2^*(y)}
F\left(f_1^*(x), s\right)\, ds
\right)j\left(|x-y|\right)\, dx\, dy
\end{eqnarray*}
where $\Psi(s_1,s_2)=\int_0^{s_2}F\left(s_1, u\right)\, du$. 

Applying Tonelli's theorem (see (3.0)), we obtain: 
\begin{equation*}
J\left(f_1^*, f_2^*\right)=\int_0^\infty\int_{\mathbb{R}^n}\int_{\mathbb{R}^n}j\left(|x-y|\right)\indicator_{\{y\in\mathbb{R}^n: f_2^*(y)\geq s\}}F\left(f_1^*(x),s\right)\, dy\, dx\, ds.
\end{equation*}

Setting $u(x, s)=\int_{\mathbb{R}^n}\indicator_{\{y\in\mathbb{R}^n: f_2^*(y)\geq s\}}j\left(|x-y|\right)\, dy$, it follows from lemma~\ref{lem21} that $u$ is radial and radially decreasing with respect to $x$ for every fixed $s$.  
\begin{equation*}
J\left(f_1^*, f_2^*\right)=\int_0^\infty\int_{\mathbb{R}^n}u(x, s)F\left(f_1^*(x),s\right)\, dx\, ds.
\end{equation*}

Now for a fixed $s_1\geq 0$, $\Psi(s_1, x_2)-\Psi(s_1, x_1)=\int_{x_1}^{x_2}F(s_1, t)\, dt\geq 0$ for $x_2\geq x_1$; from which we deduce that $F(s_1, t)$ is nonnegative for almost every $t\geq 0$. (3.0)

On the other hand, $0=\Psi(0, s_2)=\int_0^{s_2}F(0, u)\, du$. By the nonnegativity of $F$, we conclude that $F(0, s)=0$ for almost every $s\geq 0$. 

Moreover ($\Psi$3) says that: $\Psi(tx, d)-t\Psi(x, d)-\Psi(tx, c)+t\Psi(x, c)\leq 0$ for every $x\geq 0$, $d\geq c \geq 0$.  

Integrating this inequality, we have $\int_c^d F(tx, u)-t F(x, u)\, du\geq 0$ for every $x\geq 0$; $d\geq c\geq 0$.

Hence $F(tx, u)\leq t F(x, u)$ for all $x\geq 0$, $t\in]0, 1[$ and almost every $u\geq 0$. 

This shows that for almost every $s\geq 0$, the function $u(x, s)F\left(f_1^*(x), s\right)$ satisfies all the hypotheses of lemma~\ref{lem23}, consequently: 

For almost every $s\geq 0$ $\int_{\mathbb{R}^n}u(x, s)F\left(f_1^*(x),s\right)\, dx \leq\int_{\mathbb{R}^n}u(x, s)F\left(k_1\indicator_{B_1}(x),s\right)\, dx$ and 
\begin{equation}\label{eq31}
J\left(f_1^*, f_2^*\right)\leq J\left(k_1\indicator_{B_1}, f_2^*\right). 
\end{equation}

Using the same argument, we easily conclude that
\begin{equation}\label{eq32}
J\left(k_1\indicator_{B_1}, f_2^*\right)\leq J\left(k_1\indicator_{B_1}, k_2\indicator_{B_2}\right). 
\end{equation}

By \cite[Theorem 2]{c4} we know that: 
\begin{equation}\label{eq33}
J\left(f_1, f_2\right)\leq J\left(f_1^*, f_2^*\right). 
\end{equation}

Combining these three inequalities, we obtain: 
\begin{equation*}
J\left(f_1, f_2\right)\leq J\left(f_1^*, f_2^*\right)\leq J\left(k_1\indicator_{B_1}, f_2^*\right)\leq J\left(k_1\indicator_{B_1}, k_2\indicator_{B_2}\right). 
\end{equation*}

If in addition, we have strict inequality in ($\Psi$2) and ($\Psi$3), $j$ is strictly decreasing and $J\left(f_1, f_2\right)<\infty$ for any $f_1, f_2\in C$ then \cite[Theorem 2]{c4} asserts that equality occurs in~(\ref{eq33}) if and only if there exists $x_0\in\mathbb{R}^n$ such that $f_1=f_1^*(\cdot-x_0)$ and $f_2=f_2^*(\cdot-x_0)$. 

On the other hand, by lemma~\ref{lem23}, equality occurs in~(\ref{eq31}) if and only if $f_1^*=k_1\indicator_{B_1}$.  Similarly equality holds in~(\ref{eq32}) if and only if $f_2^*=k_2\indicator_{B_2}$. 

Conclusion: we have proved that for any absolutely continuous function $\Psi$ satisfying ($\Psi$1), ($\Psi$2), ($\Psi$3) with a kernel $j$ satisfying (j1) $\left(k_1\indicator_{B_1}, k_2\indicator_{B_2}\right)$ is a maximizer of $J$ under the constraint $C$.  If additionally ($\Psi$2), ($\Psi$3) hold with strict inequality $j$ is strictly decreasing and $J\left(f_1, f_2\right)<\infty$ for all $\left(f_1, f_2\right)\in C$ then $\left(k_1\indicator_{B_1}, k_2\indicator_{B_2}\right)$ is the unique maximizer of (P1) (up to a translation). 

\textbf{Remark 1:} $\Psi$ is a nondecreasing
function with respect to each variable, it is then of bounded variations. 
The absolute continuity is then equivalent to its continuity 
and the fact that it satisfies the N-Luzin property. 

\textbf{Remark 2:} We can remove condition ($\Psi$1)
from our theorem by modifying ($\Psi$3) and 
adding an integrability assumption in a same way 
as~\cite[Proposition 3.2]{c8}. 

\textbf{Part 2:} $\Psi$ is bounded.  

Applying lemma~\ref{lem22}, we know that there exist 
$\varphi_1, \varphi_2$ such that 
$\Psi(s_1, s_2)=\tilde\Psi\left(\varphi_1(s_1), \varphi_2(s_2)\right)$, 
where $\tilde\Psi$ is Lipschitz continuous with respect to each variable, there exist a function 
$\tilde{F}$ defined on $\mathbb{R}_+$ such that 
$\tilde\Psi(s_1, s_2)=\int_0^{s_2}\tilde{F}(s_1, u)\, du$. 
\begin{eqnarray*}
J\left(f_1, f_2\right)&=&\int_{\mathbb{R}^n}\int_{\mathbb{R}^n}\tilde\Psi\left(\varphi_1\left(f_1^*(x)\right), \varphi_2\left(f_2^*(x)\right)\right)j\left(|x-y|\right)\, dx\, dy\\
&=&\int_0^\infty\left(\int_{\mathbb{R}^n}\nu(x, s)\tilde{F}\left(\varphi_1\left(f_1^*(x)\right), s\right)\right)\, dx\, ds 
\end{eqnarray*}
where $\nu(x, s)=\int_{\mathbb R^n}\indicator_{\{y\in\mathbb{R}^n:\varphi_2\left(f_2^*(y)\right)\geq s\}}j\left(|x-y|\right)\, dy$. 
The function $\indicator_{\{y\in\mathbb R^n:\varphi_2\left(f_2^*(y)\right)\geq s\}}$ is Schwarz-symmetric for every $s$ since $\varphi_2$ is nondecreasing.  
We can then apply Part 1 and the result follows.  

\textbf{Remark 3:} Here we cannot obtain a uniqueness result since $\varphi_1$ and $\varphi_2$ do not inherit the strict monotonicity properties of $\Psi$. 

\textbf{Part 3:} $\Psi$ is not bounded.  

For $L>0$, set $\Psi^L\left(s_1, s_2\right)=\Psi\left(\min(s_1, L), \min(s_2, L)\right)$.  It is easy to check that $\Psi^L$ inherits all the properties of $\Psi$ stated in our result.  Moreover Part 2 applies to $\Psi^L$ since it is a bounded function.  Noticing that $\Psi^L\rightarrow\Psi$, the monotone convergence theorem enables us to conclude. 


\begin{thebibliography}{9}


\bibitem{c1}
E.~H. Lieb. 
\newblock Existence and uniqueness of the minimizing solution of the Choquard's nonlinear equation.
\newblock {\em Studies in Applied Mathematics}, 57:93--105, 1977.

\bibitem{c2}
A. Burchard. 
\newblock Cases of equality in Riesz rearrangement inequality.
\newblock {\em Ann. Math. (2)}, 143:499--527, 1996.

\bibitem{c3}
E.~H. Lieb and M.~Loss.
\newblock {\em Analysis}, volume~14 of {\em Graduate Studies in Mathematics}.
\newblock American Mathematical Society, Providence, RI, second edition, 2001.

\bibitem{c4}
A.~Burchard and H.~Hajaiej.
\newblock Rearrangement inequalities for functionals with monotone integrands.
\newblock {\em J. Funct. Anal.}, 233(2):561--582, 2006.

\bibitem{c5}
G.~R. Burton. 
\newblock Vortex-rings of prescribed impulse.
\newblock {\em Math. Proc. Cambridge Philos. Soc.}, 134(3):515--528, 2003.

\bibitem{c6}
A. Sklar. 
\newblock Functions de r\'epartition \`a $n$ dimensions et leurs marges. 
\newblock {\em Inst. Statist. Univ. Paris}, 8:229--231, 1959.

\bibitem{c7}
A. Burchard and M. Schmuckenschl\"ager. 
\newblock Comparison theorems for exist times.
\newblock {\em Geom. Funct. Anal.}, 11:651--692, 2001.

\bibitem{c8}
C. Draghici and H. Hajaiej. 
\newblock Uniqueness and characterization of maximizers of integral functionals with constraints. 
\newblock Preprint. 

\end{thebibliography}

\end{document}